\documentclass[12pt]{amsart}

\usepackage{amsmath}
\usepackage{amssymb}
\usepackage{sty1}

\makeatletter
\@addtoreset{equation}{section}
\makeatother

\def\aqg{\mathrm{A}_q(G)}
\def\coninftyg{\mathrm{Conv}^\infty(G)}
\def\contwog{\mathrm{Conv}^2(G)}
\def\fal#1{\mathrm{A}(#1)}
\def\falg{\mathrm{A}(G)}
\def\fsalg{\mathrm{B}(G)}
\def\falcg{\mathrm{A}_c(G)}
\def\falh{\mathrm{A}(H)}
\def\blone#1{\mathrm{L}^1(#1)}
\def\bloneg{\mathrm{L}^1(G)}
\def\blinfty#1{\mathrm{L}^\infty(#1)}
\def\blinftyg{\mathrm{L}^\infty(G)}
\def\blpg{\mathrm{L}^p(G)}

\def\bltwog{\mathrm{L}^2(G)}
\def\bltwo#1{\mathrm{L}^2(#1)}
\def\dist{\mathrm{dist}}
\def\idealh{\mathrm{I}(H)}
\def\lebfal#1{\mathrm{LA}(#1)}
\def\lebfalg{\mathrm{LA}(G)}
\def\sa{\mathrm{S}\mathcal{A}}
\def\sonea#1{\mathrm{S}^1\!\mathrm{A}(#1)}
\def\soneag{\mathrm{S}^1\!\mathrm{A}(G)}
\def\soneg{\mathrm{S}^1(G)}
\def\soneah{\mathrm{S}^1\!\mathrm{A}(H)}
\def\spag{\mathrm{S}^p\!\mathrm{A}(G)}
\def\spaqg{\mathrm{S}^p\!\mathrm{A}_q(G)}
\def\stwoag{\mathrm{S}^2\!\mathrm{A}(G)}
\def\vng{\mathrm{VN}(G)}
\def\vnhg{\mathrm{VN}_H(G)}
\def\vnh{\mathrm{VN}(H)}

\def\supp{\mathrm{supp}}

\def\mat#1#2{\mathrm{M}_{#1}(#2)}
\def\matm#1{\mathrm{M}_m(#1)}
\def\matn#1{\mathrm{M}_n(#1)}
\def\smat#1{\mathrm{M}_{#1}}
\def\smatn{\mathrm{M}_n}
\def\wcbop#1{\mathcal{CB}^\sigma(#1)}

\def\cbnorm#1{\left\|#1\right\|_{cb}}
\def\lonenorm#1{\left\|#1\right\|_{\mathrm{L}^1}}
\def\ltwonorm#1{\left\|#1\right\|_{\mathrm{L}^2}}
\def\linftynorm#1{\left\|#1\right\|_{\mathrm{L}^\infty}}
\def\linftynormn#1{\left\|#1\right\|_{\mathrm{M}_n(\mathrm{L}^\infty)}}
\def\fnorm#1{\left\|#1\right\|_{\mathrm{A}}}
\def\fsnorm#1{\left\|#1\right\|_{\mathrm{B}}}
\def\fnormm#1{\left\|#1\right\|_{\mathrm{M}_m(\mathrm{A})}}
\def\fnormn#1{\left\|#1\right\|_{\mathrm{M}_n(\mathrm{A})}}
\def\lfnorm#1{\left\|#1\right\|_{\mathrm{LA}}}

\def\sonenorm#1{\left\|#1\right\|_{\mathrm{S}^1}}
\def\sonefnorm#1{\left\|#1\right\|_{\mathrm{S}^1\!\mathrm{A}}}
\def\sonefdnorm#1{\left\|#1\right\|_{\mathrm{S}^1\!\mathrm{A}^*}}
\def\sonefnormn#1{\left\|#1\right\|_{\mathrm{M}_n(\mathrm{S}^1\!\mathrm{A})}}
\def\sonefdnormn#1{\left\|#1\right\|_{\mathrm{M}_n(\mathrm{S}^1\!\mathrm{A}^*)}}
\def\sonefnormmn#1{\left\|#1\right\|_{\mathrm{M}_{mn}(\mathrm{S}^1\!\mathrm{A})}}
\def\stwofnorm#1{\left\|#1\right\|_{\mathrm{S}^2\!\mathrm{A}}}
\def\vnnorm#1{\left\|#1\right\|_{\mathrm{VN}}}
\def\vnnormn#1{\left\|#1\right\|_{\mathrm{M}_n(\mathrm{VN})}}

\begin{document}

\newtheorem{sonefisosa}{Proposition}[section]

\newtheorem{symbai}{Lemma}[section]
\newtheorem{convolvers}[symbai]{Proposition}
\newtheorem{soneagdual}[symbai]{Theorem}
\newtheorem{soneagdual1}[symbai]{Corollary}

\newtheorem{distance}{Theorem}[section]
\newtheorem{restricting}[distance]{Theorem}
\newtheorem{tpfailure}[distance]{Corollary}
\newtheorem{averaging}[distance]{Theorem}
\newtheorem{averaging1}[distance]{Lemma}
\newtheorem{averaging2}[distance]{Corollary}

\newtheorem{compactgroup}{Proposition}[section]
\newtheorem{sonefamenable}[compactgroup]{Theorem}
\newtheorem{sonefamenable1}[compactgroup]{Corollary}
\newtheorem{hypertauberian}[compactgroup]{Theorem}
\newtheorem{weakamenable}[compactgroup]{Corollary}

\title[Operator Segal Algebras]
{Operator Segal Algebras in Fourier Algebras}

\author{Brian E.\ Forrest, Nico Spronk and Peter J.\ Wood}

\begin{abstract}
Let $G$ be a locally compact group, $\falg$ its Fourier algebra
and $\bloneg$ the space of Haar integrable functions on $G$.
We study the Segal algebra $\soneag=\falg\cap\bloneg$ in
$\falg$.  It admits an operator space structure which makes
it a completely contractive Banach algebra.  We compute
the dual space of $\soneag$.  We use it show that 
restriction operator $u\mapsto u|_H:\soneag\to\falh$, for
some non-open closed subgroups $H$, is a surjective complete qutient map.  
We also show that if $N$ is a non-compact
closed subgroup, then the averaging operator $\tau_N:\soneag\to\blone{G/N}$,
$\tau_Nu(sN)=\int_N u(sn)dn$ is a surjective complete quotient map.
This puts an operator space perspective on the philosophy that
$\soneag$ is ``locally $\falg$ while globally $\mathrm{L}^1$''.
Also, using the operator space structure
we can show that $\soneag$ is operator amenable exactly
when when $G$ is compact; and we can show that it is always operator
weakly amenable.  To obtain the latter fact, we use E.\ Samei's
theory of hyper-Tauberian Banach algebras.  
\end{abstract}

\maketitle

\footnote{{\it Date}: \today.


Research of the first named author supported by NSERC Grant 90749-00.
Research of the second named author supported by NSERC Grant 312515-05.
Research of the third named author supported by NSERC Grant 249726-02.}


\section{Operator Segal Algebras}

\subsection{Notation}
For any Banach space $\fX$ we let $\fB(\fX)$ denote the Banach algebra
of bounded linear operators from $\fX$ to itself, and $\ball{\fX}$
the set of all vectors of norm not exceeding $1$.

For details on classical harmonic analysis, we use \cite{hewittrI,reiter}.  
We will always let
$G$ denote a locally compact group with a fixed left invariant
Haar measure $m$.  For $1\leq p\leq \infty$,
$\mathrm{L}^p(G)$ is the usual L$^p$-space with respect to $m$.
If $f,g$ are Borel measurable functions
and $s\in G$, we have for almost every $t\iin G$ that
\[
s\con f(t)=f(s^{-1}t),\quad f\con g(t)=\int_G f(s)s\con g(t)dt\; \aand\;
\check{f}(t)=f(t^{-1})
\]
denotes the left group action, convolution (when the integrand makes sense)
and inversion.  
We note that $\bloneg$ is a Banach algebra with respect to convolution.

We let $\falg$ and $\fsalg$ denote the {\it Fourier} and 
{\it Fourier-Stieltjes algebras} of $G$ which are Banach algebras of continuous
functions on $G$ and were introduced in \cite{eymard}.
We recall, from that article, that $\falg$ consists exactly of
functions on $G$
of the form $u(s)=\inprod{\lam(s)f}{g}=\bar{g}\con\check{f}(s)$,
where $\lam:G\to\fB(\bltwog)$ is the left regular representation given
by $\lam(s)f=s\con f$.  The dual of $\falg$ is the von Neumann
algebra $\vng$, which is generated by $\lam(G)$ in $\fB(\bltwog)$.

Our standard references for operator spaces are \cite{effrosrB,pisier}.
An operator space is a complex Banach space $\fV$ equipped with
an {\it operator space structure}:
for each space of matrices $\matn{\fV}$ with entries in $\fV$, $n\in\En$,
we have a norm $\norm{\cdot}_{\matn{\fV}}$, and the norms satisfy
Ruan's axioms in addition to that $\norm{\cdot}_{\mat{1}{\fV}}$ is
the norm on $\fV=\mat{1}{\fV}$.  A map $T$ from $\fV$ to another operator space
$\fW$ is said to be {\it completely bounded} if the family of linear operators
$[v_{ij}]\mapsto [Tv_{ij}]:\matn{\fV}\to\matn{\fW}$ is uniformly bounded
over $n$.  If $\fA$ is an algebra and an operator
space, for which $\fV$ is a left module over $\fA$, we say $\fV$
is a {\it completely bounded $\fA$-module} if there is $C>0$ so for each
$[a_{ij}]\iin\matn{\fA}$ and each $[v_{kl}]\iin\matm{\fV}$ we have
\[
\norm{[a_{ij}v_{kl}]}_{\mat{nm}{\fV}}\leq C
\norm{[a_{ij}]}_{\matn{\fA}}\norm{[v_{kl}]}_{\matm{\fV}}.
\]
We say $\fV$ is a {\it completely contractive $\fA$-module}
if we can set $C=1$.
This is the same as asserting that the module multiplication,
extends to a map on the {\it operator projective tensor product}
$\fA\what{\otimes}\fV\to\fV$ is bounded at all matrix levels by $C$.
We say $\fA$ is a {\it completely bounded (contractive) Banach algebra}
if it itself is a completely bounded (contractive) $\fA$-module.
We note that any C*-algebra admits a canonical operator
space structure.  
The algebras $\bloneg$ and $\falg$ will always
have the {\it standard predual} structures (see \cite{blecher}), in
their respective roles as the preduals of $\blinftyg$ and $\vng$.
With these operator space structures, these are completely contractive 
Banach algebras.

\subsection{Abstract Operator Segal Algebras}
Let $\fA$ be a completely contractive Banach algebra.
A {\it (contractive) [left] operator Segal algebra} 
in $\fA$ is a dense [left] ideal
$\sa$ equipped with an operator space structure 
$\{\norm{\cdot}_{\matn{\sa}}:\matn{\sa}\to\Ree^{\geq 0}\}$
under which

\smallskip
{\bf (OSA1)} $(\sa,\norm{\cdot}_{\sa})$ is a Banach space,

\smallskip
{\bf (OSA2)} the identity map $\sa\hookrightarrow\fA$ is completely bounded
(contractive), and

\smallskip
{\bf (OSA3)} $\sa$ is a completely bounded (contractive) $\fA$-bimodule.

\smallskip
\noindent Note that (OSA2) and (OSA3) imply that $\sa$ is 
a completely bounded Banach algebra.  Moreover $\sa$ is a completely 
contractive Banach algebra if the associated maps and module actions are
completely contractive.  

Let us see that operator Segal algebras are reasonably common. 

{\bf (i)} For any Banach space $\fX$, let $\max\!\fX$
denote the maximal operator space whose underlying Banach space
is $\fX$.  If $\fA$ is any Banach algebra with abstract
Segal subalgebra $\sa$, then $\max\sa$ is an operator Segal algebra
in $\max\!\fA$.

{\bf (ii)} Let $\fH$ be a Hilbert space, $\fK(\fH)$ the C*-algebra
of compact operators on $\fH$ and for $1\leq p<\infty$, $\fS^p(\fH)$
be the Shatten $p$-class operators. Since we have $\fS^1(\fH)\cong
\fK(\fH)^*$ under the dual pairing $(s,k)\mapsto \mathrm{trace}(sk)$,
we see that $\fS^1(\fH)$ is a completely contractive $\fK(\fH)$-module
in its dual operator space structure,
and hence is an operator Segal algebra.  For $1<p<\infty$, we 
assign to $\fS^p(\fH)$ the interpolated operator space structure
$\fS^p(\fH)=(\fK(\fH),\fS^1(\fH))_{1/p}$.  See \cite{pisier} for
more on this.  By the functorial properties
of operator interpolation we can verify that (OSA2) and (OSA3) obtain
for $\fS^p(\fH)$, and, in fact, $\fS^p(\fH)$ is a contractive
Segal algebra in $\fK(\fH)$.

{\bf (iii)} Similarly to (ii), above, we can see that if $I$ is any set and
$c_0(I)$ denotes the C*-algebra of functions on $I$ vanishing at infinity,
then $\ell^1(I)\cong c_0(I)^*$ is a contractive operator Segal algebra
in $c_0(I)$.  With the operator interpolation structure
$\ell^p(I)=(c_0(I),\ell^1(I))_{1/p}$ ($1\leq p\leq\infty$) we
have that $\ell^p(I)$ is a contractive operator Segal algebra in
$c_0(I)$.

\smallskip
\subsection{The 1-Segal Fourier algebra}
We define the {\it 1-Segal Fourier algebra} to be the space
\[
\soneag=\falg\cap\bloneg.
\]
For $u\iin\soneag$ we let
\[
\sonefnorm{u}=\fnorm{u}+\lonenorm{u}.
\]
In \cite{ghahramanil} this space is denoted $\lebfalg$.  It is shown
in Lemma 1.1 of that article that it is complete; and in Proposition
2.5 that it is a Segal algebra in $\falg$.  We will reserve the notation
$\lebfalg$ for $\soneag$, when it is treated as a Segal algebra in $\bloneg$,
and call it the {\it Lebesgue-Fourier algebra}.  If $G$ is abelian, with dual 
group $\hat{G}$, then there is an isometric algebra isomorphism
$\soneag\cong\lebfal{\hat{G}}$.

Let us assign to $\soneag$ a natural operator space structure.  The norm
on $\soneag$ was gained via the embedding $\soneag\hookrightarrow
\falg\oplus_1\bloneg:u\mapsto (u,u)$.  
The space $\falg\oplus_1\bloneg$ is the predual
of the von Neumann algebra $\vng\oplus_\infty\blinftyg$ and as such
inherits a natural operator space structure \cite{blecher,effrosrB}.
Thus we identify the matrix space $\matn{\soneag}$ as a subspace of
$\falg\oplus_1\bloneg\cong\wcbop{\vng\oplus_\infty\blinftyg,\smatn}$.
Hence if $[u_{ij}]\in\matn{\soneag}$ we obtain
\[
\sonefnormn{[u_{ij}]}=\sup\left\{
\norm{\bigl[\dpair{T_{pq}}{u_{ij}}+
\dpair{\vphi_{pq}}{u_{ij}}\bigr]}_{\smat{nr}}:
\begin{matrix} [T_{pq}]\in\ball{\mat{r}{\vng}}  \\ 
[\vphi_{pq}]\in\ball{\mat{r}{\blinftyg}} \\ r\in\En \end{matrix}\right\}
\]
where $\dpair{T_{pq}}{u_{ij}}$ indicates the $\vng$-$\falg$ dual pairing
and $\dpair{\vphi_{pq}}{u_{ij}}=\int_G\vphi_{pq}(s)u_{ij}(s)ds$ indicates
the $\blinftyg$-$\bloneg$ dual pairing, for each index quartuple $p,q,i,j$.
It follows immediately that $\sonefnormn{[u_{ij}]}\geq\fnormn{[u_{ij}]}$,
so (OSA2) is satisfied.  It remains to check the axiom (OSA3).
Suppose $[v_{kl}]\in\ball{\mat{m}{\falg}}$.  Then
$\norm{[v_{kl}]}_{\mat{m}{\mathrm{L}^\infty}}\leq\fnormm{[v_{kl}]}\leq 1$,
since the injection $\falg\hookrightarrow\blinftyg$ is a contraction,
hence a complete contraction, so 
$\norm{[\vphi_{pq}v_{kl}]}_{\mat{rm}{\mathrm{L}^\infty}}\leq 1$
for each $[\vphi_{pq}]\iin\ball{\mat{r}{\blinftyg}}$.  Hence we have that
\begin{align*}
\sonefnormmn{[v_{kl}u_{ij}]}
&=\sup\left\{\norm{\bigl[\dpair{T_{pq}}{v_{kl}u_{ij}}
+\dpair{\vphi_{pq}}{v_{kl}u_{ij}}\bigr]}_{\smat{nr}}\right\} \\
&=\sup\left\{\norm{\bigl[\dpair{T_{pq}v_{kl}}{u_{ij}}
+\dpair{\vphi_{pq}v_{kl}}{u_{ij}}\bigr]}_{\smat{nr}}\right\} \\
&\leq\sup\left\{\norm{\bigl[\dpair{T_{pq}}{u_{ij}}
+\dpair{\vphi_{pq}}{u_{ij}}\bigr]}_{\smat{nr}}\right\} \\
&=\sonefnormn{[u_{ij}]}
\end{align*}
where the suprema are taken over choices
$[T_{pq}]\iin\ball{\mat{r}{\vng}}$, $[\vphi_{pq}]$ in 
$\ball{\mat{r}{\blinftyg}}$
and varying $r\iin\En$.  Thus if $[v_{kl}]\in\ball{\mat{m}{\falg}}$ and
$[u_{ij}]\in\matn{\soneag}$ then
\[
\sonefnormmn{[v_{kl}u_{ij}]}
\leq
\fnormm{[v_{kl}]}\sonefnormn{[u_{ij}]}.
\]
Collecting these facts together we obtain

\begin{sonefisosa}\label{prop:sonefisosa}
$\soneag$ is a contractive operator Segal algebra in $\falg$.
\end{sonefisosa}

The convolution algebra $\lebfalg$ was shown in \cite{ghahramanil}
to be a left Segal algebra in $\bloneg$.  It is immediate that under
the operator space structure developed above that the injection
$\lebfalg\hookrightarrow\bloneg$ is a complete contraction
and, since $\bloneg$ has the maximal operator space structure,
that $\lebfalg$ is a completely contractive $\bloneg$-module.
Thus $\lebfalg$ is an operator Segal algebra as well.

\smallskip
\subsection{Other Segal algebras}\label{ssec:spaq}
If $1<p<\infty$ we define the {\it Segal $p$-Fourier algebra} by
\[
\spag=\falg\cap\blpg.
\]
It is standard to show that if we embed $\spag\hookrightarrow
\falg\oplus_1\blpg$, then we obtain a Segal algebra in $\falg$.
If $G$ is abelian with dual group $\hat{G}$, then $\spag$ is
the Segal algebra in $\blone{\hat{G}}$ given by
$\{f\in\blone{\hat{G}}:\hat{f}\in\blpg\}$; see \cite{reiter}. 
If we admit on $\blpg$ 
an operator space structure for which it is a completely
contractive $\blinftyg$-module -- we might refer to this as
an ``$\mathrm{L}^\infty$-homogeneous operator space structure'', 
then we may assign an operator
space structure on $\spag$ analogously to that which we
assigned to $\soneag$, above.  The result is still a contractive
operator Segal algebra.  If $p=2$, there are many candidate
operator space structures including row, column and the
``operator Hilbert space'' $\mathrm{O}\bltwog$.  For any
$p$, the interpolated structures $\blpg=
\bigl(\min\blinftyg,\max\bloneg\bigr)_{1/p}$ suffice.
See \cite{pisier} for information on interpolation
and $\mathrm{O}\bltwog$.
By \cite[Theo.\ 3.5]{lambertnr} or \cite{lambert} there are 
certain row and column operator space structures
$\mathrm{ROW}(\blpg)$ and $\mathrm{COL}(\blpg)$
which are also $\mathrm{L}^\infty$-homogeneous operator space structures.

For $0<q<\infty$ the Fig\`{a}-Talamanca-Herz algebra
$\aqg$ has be shown in \cite{lambertnr} to admit an operator
space structure under which is is a completely bounded
Banach algebra.  If $1<p<\infty$ we can define the 
{\it Segal $p,q$-Fig\`{a}-Talamanca-Herz 
algebra} by
\[
\spaqg=\aqg\cap\blpg.
\]
If on $\blpg$ we admit a $\mathrm{L}^\infty$-homogeneous operator space 
structure, then similarly as above we obtain an operator Segal algebra
in $\aqg$.  Note that $\spaqg$ may not be a completely {\it contractive}
Banach algebra as $\aqg$ is not known to be.

\section{Dual Spaces}

In this section we develop the dual of $\soneag$ and
apply this to determining some restriction and averaging theorems.

We begin with a useful lemma.
We denote the group action of right translation by
\[
t\mult f(s)=f(st).
\]
for $t,s\iin G$, where $f$ is any function on $G$.
If $\soneg$ is a left Segal algebra in $\bloneg$, 
then we say that $\soneg$ has continuous right translations if for any
$u\iin\soneg$, $t\mult u\in\soneg$, and $t\mapsto t\mult u:G\to\soneg$ is
continuous.  For example, if $u\in\lebfalg$ we have
\[
\lfnorm{t\mult u-u}=\lonenorm{t\mult u-u}+\fnorm{t\mult u-u}
\overset{t\to e}{\longrightarrow}0.
\]
However, the right action of $G$ on $\lebfalg$ is isometric (bounded)
if and only if $G$ is unimodular -- in which case we say $\lebfalg$
is a {\it symmetric} Segal algebra in $\bloneg$.  Indeed $\lfnorm{t\mult u}
=\Del(t)\lonenorm{u}+\fnorm{u}$, where $\Del$ is the Haar modular
function.

\begin{symbai}\label{lem:symbai}
Let $\fU$ denote a neigbourhood basis of relatively compact
symmetric neighbourhoods
of the identity $e$ in $G$, which is a directed set via reverse inclusion.
For each $U\iin\fU$ we let $e_U=\frac{1}{m(u)}1_U$ (normalised indicator
function).  Then if $\soneg$ is any Segal
algebra in $\bloneg$ with continuous right translations, 
then, for any $u\iin\soneg$, 
$u\con e_U\in\soneg$ for each $U$ and
$\lim_{U\in\fU}\sonenorm{u\con e_U-u}=0$.
\end{symbai}

\proof 
If $u\in\soneg$ then for each $U\iin\fU$ and
almost every $s\iin G$ we have
\[
u\con e_U(s)=\int_G u(t)e_U(t^{-1}s)dt
=\int_G u(st)e_U(t)dt=\frac{1}{m(U)}\int_U f(st)dt
\]
where we used symmetry of $U$ to obtain that $\check{e}_U=e_U$.
Since right translation is continuous 
on $G$, and $U$ is relatively compact, we may regard
\[
u\con e_U=\frac{1}{m(U)}\int_U t\mult u\,dt
\]
as a Bochner integral, converging in $\soneg$.  We then obtain
\begin{align*}
\sonenorm{u\con e_U-u} &=
\sonenorm{\frac{1}{m(U)}\int_U(t\mult u-u)dt} \\
&\leq\frac{1}{m(U)}\int_U\sonenorm{t\mult u-u}dt 
\leq\sup_{t\in U}\sonenorm{t\mult u-u}\overset{U\in\fU}{\longrightarrow}
0.
\end{align*}
We note that $\{e_U\}_{U\in\fU}$ is a well-known symmetric bounded approximate
identity in $\bloneg$.  \endpf

\subsection{$\mathrm{L}^\infty$-convolvers} Since the metrical structure 
on $\soneag$ is determined by an embedding $\soneag\hookrightarrow
\falg\oplus_1\bloneg$, the dual $\soneag^*$ is a quotient of
$\vng\oplus_\infty\blinftyg$ by the annihilator $\soneag^\perp$.   
The $\mathrm{L}^\infty$-convolvers allow us to describe this
annihilator.

\begin{convolvers}\label{prop:convolvers}
Let $\vphi\in\blinftyg$.  Then the following statements are equivalent:

{\bf (i)} $\displaystyle\sup\bigl\{\ltwonorm{\vphi\con f}:
f\in\bltwog\cap\bloneg^\vee\aand\ltwonorm{f}\leq 1\bigr\}<+\infty$, and

{\bf (ii)} $\displaystyle\sup\left\{{\left|\int_G \vphi(s)u(s)ds\right|}:
u\in\soneag\aand\fnorm{u}\leq 1\right\}<+\infty$.

\noindent In this case the operator $f\mapsto\vphi\con f:
\bltwog\cap\bloneg^\vee\to\bltwog$ extends uniquely to
a bounded linear operator $\Lam(\vphi)$ on $\bltwog$.  Furthermore
$\Lam(\vphi)\in\vng$, and the quantities in (i) and (ii)
are each equal to $\vnnorm{\Lam(\vphi)}$.  
\end{convolvers}

In (i), $\bloneg^\vee=\{\check{f}:f\in\bloneg\}$; which is $\bloneg$
itself exactly when $G$ is unimodular.  
By \cite[20.16]{hewittrI}, if $\vphi\in\blinftyg$ and
$f\in\bloneg^\vee$, then $\vphi\con f$ makes sense; though, in general,
it is not square integrable, even if we further assume that $f$ is.
For example, if $G$ is non-compact then
$1\con f=\left(\int_G\check{f}\,dm\right)1\not\in\bltwog$
if $\int_G\check{f}\,dm\not=0$.

We call such a $\vphi$, which satisfies the conditions of the proposition,
an {\it $\mathrm{L}^\infty$-convolver} and $\Lam(\vphi)$
its {\it convolution operator}.  We denote
\[
\coninftyg=\{\vphi\in\blinftyg:\vphi\text{ is a convolver}\}.
\]

\proof (i)$\rif$(ii) 
Since $\bltwog\cap\bloneg^\vee$ defines a dense subspace of $\bltwog$,
the operator $f\mapsto\vphi\con f:
\bltwog\cap\bloneg^\vee\to\bltwog$ is a continuous linear operator
and hence extends uniquely to a bounded linear operator
$\Lam(\vphi)\iin\bdop{\bltwog}$. 

If $\rho:G\to\bdop{\bltwog}$ is the right regular representation
given by $\rho(t)f=\Del(t)^{-1/2}t\mult f$ for $f\iin\bltwog$,
and $t\iin G$, then it is easy to check that
$(\Lam(\vphi)\rho(t)-\rho(t)\Lam(\vphi))f=0$ for each
$f\iin\bltwog\cap\bloneg^\vee$, and hence for each $f\iin\bltwog$.
This implies that $\Lam(\vphi)\in\vng$.

If $f\in\bltwog\cap\bloneg^\vee$ and $g\in\bltwog\cap\bloneg$, then
we have that the associated coefficient function satisfies
\begin{equation}\label{eq:gconf}
\inprod{\lam(\cdot)f}{g}=\bar{g}\con \check{f}\in\soneag.
\end{equation}
Moreover, in the $\vng$-$\falg$ dual pairing we obtain
\begin{align}\label{eq:convolverform}
\bigl\langle \Lam(\vphi),\bar{g}\con \check{f}\bigr\rangle
&=\inprod{\vphi\con f}{g}
=\int_G \left(\int_G \vphi(s)f(s^{-1}t)ds\right)\wbar{g(t)}dt \notag \\ 
&=\int_G \vphi(s) \left(\int_G f(s^{-1}t)\wbar{g(t)}dt\right)ds
=\int_G \vphi(s)\bar{g}\con \check{f}(s)ds
\end{align}
where the use of Fubini's Theorem is justified by our choices for
$f$ and $g$.
Now if $u\in\soneag$ then we have that $\bar{u}\in\blinftyg\cap\bloneg
\subseteq\bltwog$, so $\bar{u}\in\bltwog\cap\bloneg$.
Moreover, if $(e_U)_{U\in\fU}$ is the bounded approximate identity from
Lemma \ref{lem:symbai}, then each $e_U\in\bltwog\cap\bloneg^\vee$
with $\check{e}_U=e_U$.  Also, $\lim_{U\in\fU}\sonefnorm{u\con e_U-u}=0$.
We then have that $\inprod{\lam(\cdot)e_U}{\bar{u}}=u\con e_U$ and
\begin{equation}\label{eq:coincidence}
\dpair{\Lam(\vphi)}{u}=\lim_{U\in\fU}\dpair{\Lam(\vphi)}{u\con e_U}
=\lim_{U\in\fU}\int_G \vphi(s)u\con e_U(s)ds
=\int_G \vphi(s)u(s)ds.
\end{equation}
Hence we have that
\begin{align*}
\sup&\left\{{\left|\int_G \vphi(s)u(s)ds\right|}:
u\in\soneag\aand\fnorm{u}\leq 1\right\} \\
&\leq\sup\bigl\{|\dpair{\Lam(\vphi)}{u}|:u\in\ball{\falg}\bigr\}
=\vnnorm{\Lam(\vphi)}<+\infty.
\end{align*}

(ii)$\rif$(i) Since $\bltwog\cap\bloneg$ is dense in $\bltwog$ we have
for each $f\iin\bltwog\cap\bloneg^\vee$
\[
\ltwonorm{\vphi\con f}=\sup\bigl\{|\inprod{\vphi\con f}{g}|:
g\in \bltwog\cap\bloneg\aand\ltwonorm{g}\leq 1\bigr\}.
\]
Restating (\ref{eq:convolverform}) and (\ref{eq:gconf})
we obtain for
$f\in\bltwog\cap\bloneg^\vee$ and $g\in\bltwog\cap\bloneg$ that
\[
\inprod{\vphi\con f}{g}=\int_G \vphi(s)\inprod{\lam(s)f}{g}
\aand \inprod{\lam(\cdot)f}{g}\in\soneag.
\]
Thus we have that
\begin{align*}
\sup&\bigl\{\ltwonorm{\vphi\con f}:
f\in\bltwog\cap\bloneg^\vee\aand\ltwonorm{f}\leq 1\bigr\} \\
&=\sup\left\{\left|\int_G \vphi(s)\inprod{\lam(s)f}{g}ds\right|:
\begin{matrix} f\in\bltwog\cap\bloneg^\vee \\
g\in\bltwog\cap\bloneg \\
\aand \ltwonorm{f}\ltwonorm{g}\leq 1 \end{matrix} \right\} \\
&\leq\sup\left\{\left|\int_G \vphi(s)u(s)ds\right|:
u\in\soneag\aand\fnorm{u}\leq 1\right\}<+\infty.
\end{align*}

Unfortunately, we cannot determine if the formula (\ref{eq:coincidence})
obtains for any $u\iin\falg$, since if $u\not\in\soneag$, then
the limit of the integrals in cannot be expected
to hold.  \endpf

We note that if $G$ is discrete, then it is well-known that every
element of $\vng$ is an $\ell^2$-convolver, and in particular
an $\ell^\infty$-convolver, i.e.\ $\vng=\Lam(\coninftyg)$.  
This corresponds to the fact that
$\soneag=\falg\cap\ell^1(G)=\ell^1(G)$.

If $G$ is compact, then $\blinftyg\subseteq\bloneg$, and hence
$\blinftyg=\coninftyg$.  This corresponds to the fact that
$\soneag=\falg$, in this case.  Here we see that $\Lam(\coninftyg)$
is a dense $*$-subalgebra of the (reduced)
group C*-algebra $\mathrm{C}_r^*(G)$.

There are natural questions concerning the extent of $\coninftyg$.
We note that in both of the above cases that
$\coninftyg\subseteq\bltwog$.  
Moreover, we have that $\blinftyg\cap\bloneg$, which is
necessarily always contained in $\coninftyg$,
is always a subset of $\bltwog$.  Thus we ask:
{\it for which $G$ is $\coninftyg\subseteq\bltwog$?}
It is also natural to wonder about the norm closure
$\wbar{\Lam(\coninftyg)}$.  {\it When is
$\wbar{\Lam(\coninftyg)}$ a $*$-subalgebra of $\vng$?}
{\it How big is $\wbar{\Lam(\coninftyg)}$?}
We note that $\wbar{\Lam(\coninftyg)}$ always
contains the reduced C*-algebra $\mathrm{C}_r^*(G)$ of $G$.

\smallskip
\subsection{The dual of the 1-Segal Fourier algebra}

\begin{soneagdual}\label{theo:soneagdual}
There is an isometric isomorphism
\[\soneag^*\cong\vng\oplus_\infty\blinftyg/
\{(\Lam(\vphi),-\vphi):\vphi\in\coninftyg\}.
\]
In particular, $\coninftyg$ is linearly isomorphic to a norm-closed
subspace of $\vng\oplus_\infty\blinftyg$.
\end{soneagdual}

\proof We have an isometry $\soneag\hookrightarrow
\falg\oplus_1\bloneg$.  Hence by the Hahn-Banach Theorem
every continuous linear functional may be realised
as one from $\vng\oplus_\infty\blinftyg\cong(\falg\oplus_1\bloneg)^*$
via the form
\[
\bigl\langle(T,\vphi),u\bigr\rangle
=\dpair{T}{u}+\int_G \vphi u \,dm
\]
for $(T,\vphi)\iin\vng\oplus_\infty\blinftyg$ and $u\iin\soneag$,
where $\dpair{T}{u}$ denotes the $\vng$-$\falg$ dual pairing.
It follows that the annihilator of $\soneag$ will consist
exactly of those pairs $(T,\vphi)$ for which
$\dpair{T}{u}=-\int_G \vphi u\, dm$ for each $u\iin\soneag$.
But then we have  that
\begin{align*}
\sup&\left\{\left|\int_G \vphi u \,dm\right|:u\in\soneag\aand
\fnorm{u}\leq 1\right\} \\
&\leq\sup\bigl\{|\dpair{T}{u}|:u\in\ball{\falg}\bigr\}
=\vnnorm{T}<+\infty.
\end{align*}
Hence $\vphi\in\coninftyg$.  Moreover, since $\soneag$ is dense
in $\falg$, it follows that $T=-\Lam(\vphi)$. \endpf

We have found no obvious direct computation which allows us to
obtain the following corollary, in general.  We note that, by 
\cite[Theo.\ 2.1]{burnham}, $\soneag$, being a Segal algebra in $\falg$,
has spectrum $G$.  If $G$ were assumed to have an approximate unit (maybe 
unbounded!), then the corollary would be a consequence of 
\cite[Theo.\ 1.1]{burnham}.  We let $\falcg$ denote the subalgebra of
compactly supported elements of $\falg$ and recall that $\falg$
is {\it Tauberian}, as $\falcg$ is dense in $\falg$.

\begin{soneagdual1}\label{cor:soneagdual1}
$\soneag$ is an essential $\falg$-module, i.e., 
$\falg\mult\soneag$ is a dense subspace of $\soneag$.  In particular
$\soneag$ is a Tauberian Banach algebra.
\end{soneagdual1}

\proof Let $(T,-\vphi)\iin\vng\oplus_\infty\blinftyg$
be such that
\[
\dpair{T}{uv}-\int_G \vphi uv\,dm=0
\]
for every $u\iin\falg$ and $v\iin\soneag$.  Then we see from the theorem
above that $Tu=\Lam(\vphi u)$ for every $u\iin\falg$.  
It is clear that $\falcg\subset\soneag$.  Thus, for any $v\in\falcg$
we then see that 
\[
\dpair{T}{uv}=\dpair{Tu}{v}=\int_G \vphi uv\,dm
=\dpair{\Lam(\vphi)}{uv}
\]
for all $u\iin\falg$, from which it follows that 
$\dpair{T}{w}=\dpair{\Lam(\vphi)}{w}$ for any $w\in\falcg$.
Since $\falcg$ is dense in $\falg$, $T=\Lam(\vphi)$.
It follows, again for the theorem above, that $(T,-\vphi)$
is the zero functional on $\soneag$.  The Hahn-Banach Theorem
then implies that $\falg\mult\soneag$ is a dense in $\soneag$.

We observe that $\falcg=\falcg\mult\soneag$ consists of the compactly
supported elements in $\soneag$ and is dense in there too.  \endpf

\subsection{The dual of the $2$-Segal Fourier algebra}
We note that the computations in this section work analogously for 
$\stwoag$.  Let us briefly indicate how.  Since
translations are continuous on $\stwoag$ we can obtain
that $\lim_{U\in\fU}\stwofnorm{u\con e_U-u}=0$ just as
in Lemma \ref{lem:symbai}.  Then we find for a fixed
$h\iin\bltwog$ that the following two quantities

{\bf (i)} $\displaystyle\sup\bigl\{\ltwonorm{h\con f}:
f\in\bltwog\cap\bltwog^\vee\aand\ltwonorm{f}\leq 1\bigr\}$

{\bf (ii)} $\displaystyle\sup\left\{\left|\int_G hu\,dm\right|:
u\in\stwoag\aand\fnorm{u}\leq 1\right\}$

\noindent are both equal, in particular both finite if one of them is.
We call such functions $h$ $\mathrm{L}^2$-convolvers and denote the
set of them by $\contwog$.  If $h\in\contwog$ then it defines
an bounded operator $\Lam_2(h)$ on $\bltwog$ which is an element
of $\vng$.  We thus obtain the linear isometric identification
\begin{equation}\label{eq:stwoagdual}
\stwoag^*\cong\vng\oplus_\infty\bltwog/
\{(\Lam_2(h),-h):h\in\contwog\}.
\end{equation}
We note that if $G$ is an abelian group with dual group $\hat{G}$ 
and $U:\bltwog\to\bltwo{\hat{G}}$ is the Plancherel unitary,
then $U(\contwog)=\blinfty{\hat{G}}\cap\bltwo{\hat{G}}$.
{\it Is there an analogous description for $U\Lam(\coninftyg)U^*$?}

\section{Restriction and Averaging Operations}

Let us first note that if $H$ is an open subgroup of $G$, then
we have restrictions $\falg|_H=\falh$, and $\bloneg|_H=\blone{H}$.
Both restriction operations are quotient maps, the latter so provided 
the Haar measure on $H$ is the restricted Haar measure from $G$.
Hence we obtain the restriction
\[
\soneag|_H=\soneah.
\]
We shall see that for a non-open closed subgroup, this result
can be much different.

\subsection{A distance formula}
If $H$ is any closed subgroup of $G$, the restriction map
$u\mapsto u|_H:\falg\to\falh$ is also a surjective quotient
map by \cite{herz,takesakit}.  The kernel of this map is
the ideal $\idealh=\{u\in\falg:u|_H=0\}$.  Dual to this is the fact
that the subalgebra $\vnhg$, which is the weak operator 
closure of the span of $\{\lam(s):s\in H\}$ in $\vng$, is the annihilator
of $\idealh$ in $G$ and hence $*$-isomorphic to $\falh^*\cong\vnh$.

If $H$ is a closed subgroup of $G$ we say $H$ admits a
{\it bounded approximate indicator in $G$} (after \cite{aristovrs})
if there is abounded net $(u_\alp)$ from $\fsalg$ such that

{\bf (i)} $\lim_\alp u_\alp|_Hv=v$ for each $v\iin\falh$, and

{\bf (ii)} $\lim_\alp u_\alp u=0$ for each $u\in\idealh$.

\noindent By \cite[Theo.\ 3.7]{aristovrs} we can suppose that
$\fsnorm{u_\alp}\leq 1$ for each $\alp$.  We will have a bounded
approximate indicator for $H$ in $G$ provided:
\begin{itemize}
\item $H$ is neutral in $G$, i.e.\ there is a neighbourhood
basis $\fV$ of $e$ for which $VH=HV$ for each $V\iin\fV$ 
(see \cite[Prop.\ 2.2]{kanuithl}) --
this is always true if $G$ is a small invariant neighbourhood group; or
\item $G$ is amenable (see \cite[Theo.\ 1.3]{forrestkls} and 
\cite[Prop.\ 4.1]{aristovrs}) .
\end{itemize}

Since there is a completely contractive injection $\soneag\hookrightarrow
\falg$ with dense range, there is a completely contractive 
injection $\vng\hookrightarrow\soneag^*$.  Using Theorem \ref{theo:soneagdual}
we can obtain a lower bound on the norms of the range of this map:
if $T\in\vng$ then 
\begin{align*}
\sonefdnorm{T}&
=\inf\bigl\{\max\{\vnnorm{T+\Lam(\vphi)},\linftynorm{\vphi}\}:
\vphi\in\coninftyg\bigr\}  \\
&\geq\inf\bigl\{\vnnormn{T+\Lam(\vphi)}:\vphi\in\coninftyg\bigr\} \\
&=\dist_{\mathrm{VN}}\bigl(T,\Lam(\coninftyg)\bigr). 
\end{align*}
We obtain a similar bound for matricial norms:  if 
$[T_{ij}]\in\matn{\vng}$, then
\begin{equation}\label{eq:distlbd}
\sonefdnormn{[T_{ij}]}\geq
\dist_{\mathrm{M}_n(\mathrm{VN})}\bigl([T_{ij}],\matn{\Lam(\coninftyg)}\bigr).
\end{equation}
Under certain assumptions, the quantity on the right is as big as it can be.

\begin{distance}\label{theo:distance}
If $H$ is a non-open closed subgroup of $G$ which admits a contractive
approximate indicator $(u_\alp)$, then for any $[T_{ij}]\iin\matn{\vnhg}$
\[
\dist_{\matn{\mathrm{VN}}}\bigl([T_{ij}],\matn{\Lam(\coninftyg)}\bigr)
=\vnnormn{[T_{ij}]}.
\]
\end{distance}

\proof  Let us first suppose that $n=1$.  Given any $\eps>0$ find
$u\in\ball{\falg}$ for which $\supp(u)$ is compact and
\[
|\dpair{T}{u}|>\vnnorm{T}-\eps.
\]
Now let $\vphi\in\coninftyg$.  We then have for any $\alp$ that
\begin{equation}\label{eq:tplusphialp}
\dpair{T+\Lam(\vphi)}{uu_\alp}=\dpair{T}{uu_\alp}+\int_G\vphi uu_\alp\, dm.
\end{equation}
If we can find $\alp_0$ for which
\begin{equation}\label{eq:conditions}
|\dpair{T}{u}-\dpair{T}{uu_{\alp_0}}|<\eps\quad
\aand\quad
\left|\int _G\vphi uu_{\alp_0}\,dm\right|<\eps
\end{equation}
then it follows that 
\begin{equation}\label{eq:lowbd}
|\dpair{T+\Lam(\vphi)}{uu_{\alp_0}}|>\vnnorm{T}-3\eps.
\end{equation}
But then $\dist_{\mathrm{VN}}\bigl(T,\Lam(\coninftyg)\bigr)\geq \vnnorm{T}$,
whence equality must hold.

Let us verify the first inequality from (\ref{eq:conditions}).  
Since $T\in\vnhg$ there
is $T_H\iin\vnh$ such that $\dpair{T}{v}=\dpair{T_H}{v|_H}$ for each 
$v\iin\falg$.
Hence by condition (i) in the definition of $(u_\alp)$ we have
\[
\dpair{T}{uu_\alp}=\dpair{T_H}{(uu_\alp)|_H}\overset{\alp}{\longrightarrow}
\dpair{T_H}{u|_H}=\dpair{T}{u}
\]
and we can find $\alp_1$ for which the desired inequality is satisfied
for any $\alp_0\geq\alp_1$.
Now let us verify the second inequality from (\ref{eq:conditions}).  Since
$H$ is not open it is locally null \cite[20.17]{hewittrI}, so $m(S\cap H)=0$
where $S=\supp(u)$.  Let $C$ be any compact neighbourhood of $S\cap H$
with
\[
m(C)<\frac{\eps}{2\linftynorm{\vphi}+1}.
\]
Let $U$ be a symmetric neighbourhood of $e$ for which $(S\cap H) U^2
\subseteq C$.  Then if for $s\iin G$ we let
\[
v(s)=\frac{1}{m(U)}\inprod{\lam(s)1_U}{1_{(S\cap H) U}}
=\frac{m\bigl(sU\cap(S\cap H) U\bigr)}{m(U)}
\]
we obtain that
\[
\linftynorm{v}=1, \quad v|_{S\cap H}=1\quad\aand\quad
\supp(v)\subseteq (S\cap H) U^2\subseteq C.
\]
Now we have that
\[
\int_G \vphi uu_\alp\,dm=\int_G\vphi uu_\alp v\,dm
+\int_G \vphi uu_\alp(1-v)dm
\]
where $\linftynorm{uu_\alp v}\leq 1$ so
\[
\left|\int_G\vphi uu_\alp v\,dm\right|
\leq\int_C|\vphi|\,dm 
\leq\linftynorm{\vphi}m(C)<\frac{\eps}{2}
\]
while $\linftynorm{u_\alp u(1-v)}\leq\fnorm{u_\alp u(1-v)}
\overset{\alp}{\longrightarrow}0$ as $u(1-v)\in\idealh$, so
\[
\left|\int_G\vphi uu_\alp(1-v)dm\right|
\leq\int_S\linftynorm{\vphi}\linftynorm{uu_\alp(1-v)}dm
\overset{\alp}{\longrightarrow}0.
\]
Hence we can find $\alp_2$ so that the desired inequality holds
for $\alp_0\geq\alp_2$.  To obtain (\ref{eq:conditions}), we
now choose $\alp_0\geq\alp_1,\alp_2$.

We can now obtain the theorem for $n>1$.  First, given $\eps>0$
we can $m$ and
$[u_{kl}]\iin\ball{\mat{m}{\falg}}$ for which each $\supp(u_{ij})$
is compact and
\[
\norm{\dpair{T_{ij}}{u_{kl}}}_{\smat{mn}}>\vnnormn{[T_{ij}]}-\eps.
\]
Then, as in (\ref{eq:tplusphialp}), if $[\vphi_{ij}]\in\matn{\coninftyg}$
we have
\[
[\dpair{T_{ij}+\Lam(\vphi_{ij})}{u_{kl}u_\alp}]
=[\dpair{T_{ij}}{u_{kl}u_\alp}]+\left[\int_G \vphi_{ij}u_{kl}u_\alp dm\right].
\]
As above we can arrange to find an $\alp_0$ for which
\[
\left|\dpair{T_{ij}}{u_{kl}u_{\alp_0}}-\dpair{T_{ij}}{u_{kl}}\right|
<\frac{\eps}{n^2}\quad\aand\quad
\left|\int_G \vphi_{ij}u_{kl}u_{\alp_0}dm\right|<\frac{\eps}{n^2}
\]
then that the matricial analogue of (\ref{eq:lowbd}) is satisfied.
\endpf

\subsection{A restriction theorem}
We conjecture that the following theorem
holds for general closed non-open subgroups of $G$.

\begin{restricting}\label{theo:restricting}
If $H$ is a closed non-open subgroup of $G$ which admits a
bounded approximate indicator in $G$, then the restriction map
\[
u\mapsto u|_H:\soneag\to\falh
\]
is a surjective complete quotient map.
\end{restricting}

\proof The adjoint of the restriction map is the composition of
the injective $*$-homomorphism $J:\vnh\to\vnhg\subset\vng$
with the completely contractive injection
$\vng\hookrightarrow\soneag^*$.  It follows from (\ref{eq:distlbd})
and Theorem \ref{theo:distance} that for $[T_{ij}]\iin\matn{\vnh}$, we have
\begin{align*}
\sonefdnormn{[JT_{ij}]}
&\geq
\dist_{\matn{\mathrm{VN}}}\bigl([JT_{ij}],\matn{\Lam(\coninftyg)}\bigr) \\
&=\vnnormn{[JT_{ij}]}=\norm{[T_{ij}]}_{\matn{\mathrm{VN}_H}}.
\end{align*}
Hence $J:\vnhg\to\soneag^*$ is a complete isometry, so
the restriction map must be a surjective complete quotient map.
\endpf

We recall that for the operator projective tensor product $\what{\otimes}$
we have formulas
\[
\bloneg\what{\otimes}\bloneg\cong\blone{G\cross G}
\aand
\falg\what{\otimes}\falg\cong\fal{G\cross G}.
\]
It is thus a little surprising that the 1-Segal Fourier algebra $\soneag$ 
admits no such formula.

\begin{tpfailure}
If $G$ is a non-compact, non-discrete group which is either
amenable or admits small invariant neighbourhoods, then
\[
\soneag\what{\otimes}\soneag\not\cong\sonea{G\cross G}.
\]
\end{tpfailure}

\proof If the two spaces were isomorphic, then
the multiplication map $\soneag\what{\otimes}
\soneag\to\soneag\subset\falg$ would be isomorphic to
the restriction map $w\mapsto w|_{D}:\sonea{G\cross G}\to\fal{D}
\cong\falg$, where $D=\{(s,s):s\in G\}$ is the diagonal
subgroup.  Since our assumptions on $G$ allow $D$ a bounded approximate
indicator in $G\cross G$ by \cite[Theo.\ 2.4]{aristovrs},
the restriction map is surjective.  But
this contradicts that $\soneag\subsetneq\falg$. \endpf

\subsection{An averaging theorem}
We note an analogous result to Theorem \ref{theo:restricting}
which applies to $\lebfalg$, which we recall is $\bloneg\cap\falg$
{\it qua} Segal algebra of $\bloneg$.

Let $N$ be a normal subgroup of $G$, admitting left Haar integral
$\int_N\cdots dn$ (normalised if $N$ is compact).
It is well-known that the $N$-averaging operator $\tau_N:\bloneg\to
\blone{G/N}$ given by
\[
\tau_Nf(tN)=\int_N f(tn)dn
\]
for almost every $tN$ in $G/N$, is a contractive surjective algebra
homomorphism.  See \cite[III.4]{reiter}, for example.

Let us note that $\tau_N$ is a quotient map, but offer an alternative proof to
the standard one.  We let
\[
\blinfty{G\!:\!N}=\left\{\vphi\in\blinftyg:
\begin{matrix}\text{ for locally almost every }t\iin G, \\
\vphi(tn)=\vphi(t)\text{ for every }n\iin N\end{matrix}\right\}
\]
which is clearly a closed subspace of $\blinftyg$.
We note that if $\vphi\in\blinfty{G\!:\!N}$, then for locally almost
every $t\iin G$ and every $n\iin N$ we have $\vphi(nt)=
\vphi(t\mult t^{-1}nt)=\vphi(t)$.  If $q:G\to G/N$ is the quotient
map then $\vphi\mapsto\vphi\comp q$ is a linear isometry
from $\blinfty{G/N}$ onto $\blinfty{G\!:\!N}$.  Now let us
compute the adjoint of $\tau_N$:  if $\vphi\in\blinfty{G/N}$ and
$f\in\bloneg$, then
\begin{align*}
\dpair{\tau_N^*\vphi}{f}&=\dpair{\vphi}{\tau_Nf}
=\int_{G/N}\vphi(tN)\int_N f(tn)dn\,dtN \\
&=\int_{G/N}\int_N \vphi\comp q(tn)f(tn)dn\,dtN
=\int_G \vphi\comp q(t)f(t)dt=\dpair{\vphi\comp q}{f}
\end{align*}
by Weil's integral formula, so $\tau_N^*\vphi=\vphi\comp q$.
Thus $\tau_N^*:\blinfty{G/N}\to\blinftyg$ is an isometry, so $\tau_N$
is a quotient map.

If $N$ is a compact group, then we may let $\blone{G\!:\!N}$ be defined
similarly as $\blinfty{G\!:\!N}$, and we note it is a subalgebra of
$\bloneg$ which is isometrically algebraically isomorphic to $\blone{G/N}$.  
We may consider $\tau_N$ to have its range as $\blone{G\!:\!N}$.
It is well known that
$\tau_N\falg=\fal{G\!:\!N}=\{u\in\falg:u\text{ is constant on
cosets of }N\}\cong\fal{G/N}$, and that $\tau_Nu=u$ for $u\iin\falg$.  
Hence it follows that
\[
\tau_N\bigl(\lebfalg\bigr)=\lebfal{G\!:\!N}\cong\lebfal{G/N}
\]
where $\lebfal{G\!:\!N}=\{u\in\lebfalg:u\text{ is constant on
cosets of }N\}$.
This result does not hold if $N$ is not compact.  In the case that
$G$ is abelian, the following 
due to Krogstad, whose unpublished result is announced
in \cite{feichtinger}. Our proof is for general locally compact $G$.

\begin{averaging}\label{theo:averaging}
If $N$ is a non-compact closed normal subgroup of $G$, then
the $N$-averaging operator $\tau_N:\lebfalg\to\blone{G/N}$
is a surjective complete quotient map.
\end{averaging}

As with Theorem \ref{theo:restricting}, the proof relies on a distance
formula.

\begin{averaging1}\label{theo:averaging1}
If $N$ is a non-compact closed normal subgroup of $G$, and
$[\vphi_{ij}]\in\matn{\blinfty{G\!:\!N}}$, then
\[
\dist_{\matn{\mathrm{L}^\infty}}\bigl([\vphi_{ij}],
\matn{\coninftyg}\bigr)=\norm{[\vphi_{ij}]}_{\matn{\mathrm{L}^\infty}}.
\]
\end{averaging1}

\proof Let us begin with the ``scalar'' case.  Let $\vphi\in\blinfty{G\!:\!N}$
and $\eps>0$.  Let $f\in\bigl(\bloneg\cap\bltwog\bigr)^\vee$ be so
$\lonenorm{f}=1$ and $|\dpair{\vphi}{f}|>\linftynorm{\vphi}-\eps$. 
Then for any $\psi\in\coninftyg$ and any $n\in N$ we have
\begin{align*}
\dpair{\vphi+\psi}{n\con f}
&=\dpair{n^{-1}\con\vphi}{f}+\dpair{\psi}{n\con f} \\
&=\dpair{\vphi}{f}+\int_G\psi(t)\check{f}(t^{-1}n)dt=\psi\con\check{f}(n).
\end{align*}
By \cite[20.16]{hewittrI} we find that $h=\psi\con\check{f}$ is left
uniformly continuous.  Since $\psi\in\coninftyg$, $h\in\bltwog$.
It then follows that $h$ is a continuous function vanishing at $\infty$.
Indeed, if not, then there is an $\del>0$ and a net $(t_\alp)$ 
in $G$ such that $\lim_\alp t_\alp=\infty$ and for which
$|h(t_\alp)|>\del$ for each $\alp$.  Then, by uniform continuity, there
is a compact neighbourhood $U$ of the identity in $G$ such that
$|h(t_\alp s)|>\del/2$ for $s\in U$.  By dropping to a subnet, we may
assume that $t_\alp U\cap t_\beta U=\varnothing$ if $\alp\not=\beta$.
But then, selecting any finite collection $F$ of indices we find
\[
\int_G|h|^2dm\geq\sum_{\alp\in F}\int_{t_\alp U}|h|^2dm\geq
\frac{|F|\del^2}{4}
\]
which, since $|F|$ can be chosen arbitrarily large, 
contradicts that $h\in\bltwog$.

Thus, since $N$ is non-compact, we may find $n\iin N$ for which
\[
|\dpair{\vphi+\psi}{n\con f}|>\linftynorm{\vphi}-\eps.
\]
Hence $\dist_{\mathrm{L}^\infty}\bigl(\vphi,\coninftyg\bigr)=
\linftynorm{\vphi}$.

The general matricial case can be deduced from the scalar case, 
exactly as in the proof of Theorem \ref{theo:distance}.  \endpf

\noindent {\bf Proof of Theorem \ref{theo:averaging}.}
The adjoint of $\tau_N:\lebfalg\to\blone{G/N}$ is the composition
of the injective $*$-homomorphism $\vphi\mapsto\vphi\comp q:
\blinfty{G/N}\to\blinfty{G\!:\!N}\subset\blinftyg$
with the contractive inclusion $\blinftyg\hookrightarrow\lebfalg^*$.
We recall that $\lebfalg^*=\soneag^*$, which is described in Theorem
\ref{theo:soneagdual}.  Now if $[\vphi_{ij}]\in\matn{\blinfty{G/N}}$, then
we obtain, just as in (\ref{eq:lowbd}), and using the lemma above
\begin{align*}
\linftynormn{[\vphi_{ij}]}
&\geq\norm{[\tau_N^*\vphi_{ij}]}_{\matn{\mathrm{LA}^*}} \\
&\geq\dist_{\matn{\mathrm{L}^\infty}}\bigl([\vphi_{ij}\comp q],
\matn{\coninftyg}\bigr) \\
&=\linftynormn{[\vphi_{ij}\comp q]}=\linftynormn{[\vphi_{ij}]}.
\end{align*}
Hence $\tau_N^*:\blinfty{G/N}\to\lebfalg^*$ is a complete isometry,
which implies that $\tau_N:\lebfalg\to\blone{G/N}$ is a surjective complete
quotient map.  \endpf

It was shown in \cite{ghahramanil2} that if $G$ is a unimodular
group -- so $\lebfalg$ is a symmetric Segal algebra --
then {\it $\lebfalg$ is Arens regular if and only if $G$ is compact.}
Let us briefly note that this result holds without assuming {\it a priori}
that $G$ is unimodular.  We refer reader to the survey article 
\cite{duncanh} for details on, and functorial properties of, Arens regularity.

\begin{averaging2}
If $G$ is not unimodular, then $\lebfalg$ is not Arens regular.
\end{averaging2}

\proof Let us first note that if $G$ is admits a 
continuous homomorphism $\del:G\to\Ree$ for which $\ker\del$ is compact,
then for and $a>0$, $\del^{-1}([-a,a])$ is a compact neighbourhood
of the identity, invariant under inner automorphisms,
making $G$ an [IN]-group.  Then it is well known
(see \cite[12.1.9]{palmer}, for example) that $G$ is unimodular.

Hence, if $G$ is not unimodular, then $N=\ker\Delta$ is necessarily
a non-compact closed subgroup of $G$, with $G/N$ isomorphic to an infinite
subgroup of $\Ree$.  Hence $\blone{G/N}$ is not Arens regular and
$\tau_N:\lebfalg\to\blone{G/N}$ is a quotient homomorphism.  \endpf

\section{Cohomological Properties}

In this section we discuss amenability and operator amenability
of the 1-Segal Fourier algebra $\soneag$.  
We have opted to keep our presentation simple by focusing on
$\soneag$.  With suitable modifications to the proofs, we believe
that all of our results hold for $\spaqg$, which was defined
and given an operator space structure in Section \ref{ssec:spaq},
where $1\leq p<\infty$, $1<q<\infty$.

\subsection{Amenability}  Ruan's result, that $\falg$ is
{\it operator amenable} if and only if $G$ is amenable \cite{ruan},
is one of the most important results which justifies treating $\falg$
as an operator space.  We recall that a completely contractive
Banach algebra $\fA$ is {\it operator amenable} if every 
completely bounded derivation $D:\fA\to\fV^*$, where
$\fV^*$ is the operator dual space to completely bounded
$\fA$-bimodule, is inner.  We note that there exist compact groups,
for example $G=\mathrm{SO}(3)$, such that $\falg$ is
not amenable \cite{johnson}.

\begin{compactgroup}\label{prop:compactgroup}
$\soneag=\falg$ completely isomorphically if and only if $G$ is compact.
\end{compactgroup}

\proof By \cite[Prop.\ 2.6]{ghahramanil}, $\soneag=\falg$ if and only
if $G$ is compact.  Thus it remains to show that if $G$
is compact then the identity map $j:\falg\to\soneag$
is completely bounded.  We have for $[u_{ij}]\iin\matn{\falg}$ that
\[
\sonefnormn{[u_{ij}]}=\sonefnormn{[u_{ij}1]}\leq
\fnormn{[u_{ij}]}\sonefnorm{1}=2\fnormn{[u_{ij}]}.
\]
Thus $j$ is completely bounded with $\cbnorm{j}\leq 2$.
In fact, since $2=\norm{j(1)}$ we have that $2\leq\norm{j}
\leq\cbnorm{j}$, as well.  Since $j^{-1}:\soneag\to\falg$ is 
completely contractive we obtain that $\soneag=\falg$
completely isomorphically.  \endpf

Thus we obtain the main result of this section.

\begin{sonefamenable}\label{theo:sonefamenable}
$\soneag$ is operator amenable if and only if $G$ is compact.
\end{sonefamenable}

\proof ($\rif$) $\soneag$ is operator amenable only if it has
a bounded approximate identity \cite{ruan,johnsonm}.  This
happens if and only if $G$ is compact by \cite[Prop.\ 2.6]{ghahramanil}.

($\lif$) By Proposition \ref{prop:compactgroup} above, any
completely bounded $\soneag$-module $\fV$ is a completely
bounded $\falg$-module.  Moreover, every
completely bounded derivation $D:\soneag\to\fV^*$
induces a completely bounded derivation $D\comp j:
\falg\to\fV^*$.  By \cite{ruan} $D\comp j$ is inner, whence so
too is $D$.  \endpf

The main result of \cite{forrestr} states that $\falg$ is amenable
if and only if $G$ admits an abelian subgroup of finite index.
Combining this with \cite[Cor.\ 2.7]{ghahramanil} we obtain
the following.

\begin{sonefamenable1}\label{cor:sonefamenable1}
$\soneag$ is amenable if and only if $G$ is compact and
admits an abelian subgroup of finite index.
\end{sonefamenable1}

\subsection{Weak Amenability}  The theory of hyper-Tauberian
Banach algebras, developed by Samei \cite{samei}, extends very 
easily to Segal algebras.  We are grateful to E.\ Samei for
pointing out an error we made, applying his work, 
in an earlier draft of this article.

Let us first recall some basic definitions.
Let $\fA$ be a semi-simple abelian (completely contractive) Banach
algebra with Gelfand spectrum $X$.  Hence we identify $\fA$ as
a subspace of $\fC_0(X)$.  We define for $a\iin\fA$ its support
by $\supp(a)=\wbar{\{x\in X:a(x)\not=0\}}$.
We say that $\fA$ is {\it Tauberian} if
the subalgebra of compactly supported elements, $\fA_c$,
is dense in $\fA$.  
If $\fV$ is a symmetric Banach $\fA$-module
and $v\in\fV$, we define the support of $v$ over $\fA$ by
\[
\supp_\fA(v)=\{x\in X:a(x)=0\text{ whenever }a\mult v=0\}.
\]
If $\fV=\fA^*$, with the usual dual module action: $a\mult f(b)=f(ba)$,
this agrees with the usual notion of support of a linear functional.
If $\fV=\fA$, then $\supp_\fA(a)=\supp(a)$ for any $a\iin\fA$.
If $\fV$ and $\fW$ are (completely bounded) Banach $\fA$-modules
a linear operator $T:\fV\to\fW$ is called {\it $\fA$-local} if
\[
\supp_\fA(Tv)\subseteq\supp_\fA(v)\text{ for any }v\iin\fV.
\]
We say that $\fA$ is {\it (operator) hyper-Tauberian} if
every (completely) bounded local operator $T:\fA\to\fA^*$
is an $\fA$-module map.

Now suppose that the semi-simple abelian (completely contractive)
Banach algebra $\fA$ has an abstract (operator) Segal
algebra $\sa$.  By \cite[Theo.\ 2.1]{burnham}, $\sa$ is also semi-simple with 
Gelfand spectrum $X$.  

\begin{hypertauberian}\label{theo:hypertauberian}
If $\fA$ is (operator) hyper-Tauberian 
and $\sa$ is an essential $\fA$ module, i.e., $\fA\mult\sa$ is
dense in $\sa$, then $\sa$ is (operator) hyper-Tauberian.
\end{hypertauberian}

Note that the converse follows from \cite[Theo.\ 4.6]{samei},
and (OSA2).

\medskip
\proof  It is immediate that for $v\iin\sa$
\begin{equation}\label{eq:support}
\supp_\fA(v)=\supp(v)=\supp_{\sa}(v).
\end{equation}
Also, since $\sa$ is a subalgebra of $\fA$, we have that if $f\in\sa^*$
then
\begin{equation}\label{eq:dualsupport}
\supp_\fA(f)\subseteq\supp_{\sa}(f).
\end{equation}
Now suppose $T:\sa\to\sa^*$ is a (completely) bounded $\sa$-local
operator.  Then we have for $u\iin\sa$, combining (\ref{eq:dualsupport})
and (\ref{eq:support}),
\[
\supp_\fA(Tu)\subseteq\supp_{\sa}(Tu)\subseteq\supp(u)=\supp_\fA(u).
\]
Hence $T$ is also an $\fA$-local operator.  Hence by \cite[Prop.\ 2.3]{samei}
$T$ is an $\fA$-module map.  Thus it is an $\sa$-module map.  \endpf

One of the main motivations for studying hyper-Tauberian algebras
is the result \cite[Theo.\ 3.2]{samei}:  {\it if $\fA$ is (operator)
hyper-Tauberian, then it is (operator) weakly amenable.}  
We recall that $\fA$ is (operator)
weakly amenable if every (completely) bounded derivation 
$D:\fA\to\fA^*$ is inner \cite{badecd}.
We also recall that the subalgebra $\falcg$ of the
Fourier algebra of compactly supported elements lies
within $\soneag$, and is exactly the subalgebra of
compactly supported elements there. 

\begin{weakamenable}\label{cor:weakamenable}
{\bf (i)} $\soneag$ is always operator weakly amenable.

{\bf (ii)} $\soneag$ is weakly amenable if the connected component
$G_e$ of $G$ is abelian.
\end{weakamenable}

\proof In \cite[Theo.\ 7.4]{samei} it is shown that $\falg$ is always
operator hyper-Tauberian.  Thus (i) follows from \cite[Theo.\ 3.2]{samei},
whose statement was mentioned above, the preceeding theorem,
and Corollary \ref{cor:soneagdual1}.
In \cite[Theo.\ 6.5]{samei} it is shown that $\falg$ is hyper-Tauberian
if $G_e$ is abelian.  Hence (ii) follows, similarly to (i) above. 
\endpf

We note that weak amenability of $\falg$ is discussed in 
\cite{johnson} and \cite{forrestr},
and operator weak amenability in \cite{forrestw}, \cite{spronk}
and \cite{samei0}.

{
\bibliography{segalfourierbib}
\bibliographystyle{plain}
}

\smallskip
{\sc Department of Pure Mathematics, University of Waterloo,
Waterloo, ON\quad N2L 3G1, Canada}

E-mail addresses: {\tt beforres@uwaterloo.ca}, {\tt nspronk@uwaterloo.ca},
\linebreak {\tt pwood@uwaterloo.ca}

\end{document}